\documentclass[11pt]{amsart}
\usepackage{amssymb,latexsym,comment,url}
\usepackage{float}
\usepackage{tikz-cd}
\usepackage{hyperref}
\hypersetup{pdftitle={Rational Preperiodic Points of Quadratic Rational Maps over Q with Nonabelian Automorphism Groups},
  pdfauthor={Hasan Bilgili, Mohammad Sadek}}
\usepackage{caption}

\newcommand{\Char}{\operatorname{char}}

\newcommand{\PGL}{\operatorname{PGL}}

\newcommand{\Aut}{\operatorname{Aut}}

\newcommand{\isom}{ \cong }

\newcommand{\Preper}{\operatorname{PrePer}}

\newcommand{\PP}{{\mathbb P}}
\newcommand{\Q}{{\mathbb Q}}

\newcommand{\Magma}{{\sf Magma}}

\newcommand{\Mathematica}{{\sf Mathematica}}

\newenvironment{Proof}{\par\noindent{\sc Proof:}}%
                      {\hspace*{\fill}\nobreak$\Box$\par\medskip}
\newenvironment{ProofOf}[1]{\par\noindent{\sc Proof of #1:}}%
                       {\hspace*{\fill}\nobreak$\Box$\par\medskip}

\newtheorem{Proposition}{Proposition}[section]
\newtheorem{Theorem}[Proposition]{Theorem}
\newtheorem{Lemma}[Proposition]{Lemma}
\newtheorem{Corollary}[Proposition]{Corollary}
\newtheorem{Example}[Proposition]{Example}

\theoremstyle{definition}
\newtheorem{Definition}[Proposition]{Definition}
\newtheorem{Remark}[Proposition]{Remark}

\renewcommand{\baselinestretch}{1.1}
\renewcommand{\arraystretch}{1.3}

\begin{document}

\title[QUADRATIC MAPS WITH NONABELIAN AUTOMORPHISM GROUPS]
{Rational Preperiodic Points of Quadratic Rational Maps over \texorpdfstring{$\mathbb{Q}$}{Q} with Nonabelian Automorphism Groups}%

\author[H. Bilgili]%
{Hasan Bilgili}
\address{Faculty of Engineering and Natural Sciences, Sabanc{\i} University, Tuzla, \.{I}stanbul, 34956 T\"{u}rkiye}
\email{hasan.bilgili@sabanciuniv.edu}

\author[M. Sadek]%
{Mohammad Sadek}
\address{Faculty of Engineering and Natural Sciences,
 Sabanc{\i} University,
  Tuzla, \.{I}stanbul, 34956 T\"{u}rkiye}
\email{mohammad.sadek@sabanciuniv.edu}

\begin{abstract} 
Let $f:\PP^1\to\PP^1$ be a quadratic rational map defined over the rational field $\mathbb Q$ with nonabelian automorphism group. We prove that no such map has a $\mathbb Q$-rational periodic point with exact period $N\ge 4$. We also give an explicit parametrization of such maps that have $\mathbb Q$-rational
periodic points of exact period $1$, $2$, and $3$. In addition, we show that the number of $\mathbb Q$-rational preperiodic points of such a map $f$ cannot exceed $6$. As a result, we completely classify all portraits of $\mathbb Q$-rational preperiodic points for quadratic rational maps defined over $\mathbb Q$ with nonabelian automorphism group, showing that there are exactly $7$ such portraits, one of which is the empty portrait.
\end{abstract}

\maketitle

\section{Introduction} 
The uniform boundedness conjecture in arithmetic dynamics was proposed by Morton and Silverman, \cite{Mortonsilverman}, as an analogue of the uniform boundedness of torsion points of elliptic curves over a number field $K$.
Namely, there is a constant $C:=C(n,d,[K:\Q])$ such that the number of $K$-rational preperiodic points of any morphism $f:\mathbb{P}^n\to\mathbb{P}^n$, $n\ge 1$, of degree $d$ defined over $K$ is at most $ C$. Despite being a few decades old, the conjecture remains out of reach. 

In what follows, we briefly display the current knowledge for rational maps $f:\PP^1\to\PP^1$ of degree $2$ defined over $\Q$. A quadratic polynomial defined over $\Q$ cannot have a rational periodic point of exact period $4$, \cite{morton}; or exact period $5$, \cite{flynn}; or exact period $6$ assuming the Birch--Swinnerton-Dyer conjecture, \cite{stoll}. A conjecture of Poonen asserts that a quadratic polynomial defined over $\Q$ cannot have a rational periodic point of exact period $\ge 6$. In addition, if the latter conjecture holds, then the number of rational preperiodic points of a quadratic polynomial defined over $\Q$ cannot exceed $9$, \cite{poonen}. There are a few results on the boundedness of rational preperiodic points of parametric families of polynomials of arbitrary degree defined over $\Q$, see for example \cite{Ingram, Sadek, SadekUyar}. 

A rational map of degree $2$ has a nontrivial automorphism if and only if it is conjugate to a map in the unique normal form $  k\left(z + 1/z\right)$, where $k \in \mathbb{C}\setminus\{0\}$. The automorphism group is cyclic of order $2$, $\mathfrak{C}_2$, if $k \neq -\frac{1}{2}$, but is nonabelian of order $6$, $\mathfrak{S}_3$, if $k = -\frac{1}{2}$, \cite[Theorem 5.1]{Milnor}. A quadratic rational map defined over $\Q$ with automorphism group $\mathfrak{C}_2$ is proved not to have a rational periodic point of exact period $3$, \cite{Manes}. It is conjectured that such a map cannot have a rational periodic point of exact period $\ge 5$. Further, if the latter conjecture holds, then the number of rational preperiodic points of a quadratic rational map defined over $\Q$ with automorphism group $\mathfrak{C}_2$ cannot exceed $12$. In this article, we provide an analogous analysis of rational preperiodic points of quadratic rational maps defined over $\Q$ with automorphism group $\mathfrak{S}_3$.

A quadratic rational map defined over a field $K$, $\Char K\ne 2,3$, with automorphism group $\mathfrak{S}_3$ is $K$-conjugate to a rational map of the form 
\begin{equation*}
     \theta_{d,k}(z) = \frac{kz^2 - 2dz + dk}{z^2 - 2kz + d},
\end{equation*}
where $k \in K$, $d \in K\setminus\{0\}$ and $k^2\ne d$, \cite{Yasuf}. Using this explicit description of such maps, we study their rational preperiodic points. We give a complete parametrization of such maps with rational periodic points of exact period $1,2$ or $3$. 

In \cite{Levy}, the authors proved that for any rational map defined over a number field, there exists an absolute bound on the number of rational preperiodic points of any of its twists. As an application, they showed that the number of $\Q$-rational preperiodic points of a twist of $\phi(z) = z + 1/z$ cannot exceed $6$. More precisely, if $\phi_b(z) = z + \frac{b}{z}$, $b \in \mathbb{Q}^*$, then $\phi_b(z)$ has at most $6$ $\mathbb{Q}$-rational preperiodic points. Moreover, if  $\psi_b(z) = -\left(z + \frac{b}{z}\right)$, $b \in \mathbb{Q}^*$, then $\psi_b(z)$ has either $2$ or $4$ $\Q$-rational preperiodic points. A quadratic rational map with nonabelian automorphism is a twist of the map $1/z^2$, \cite{Yasuf}. Such a map is isomorphic to $1/z^2$ over a field extension of degree at most $6$. Since the preperiodic points of $1/z^2$ consist of $0, \infty$, and roots of unity, the framework of \cite{Levy} yields the upper bound $20$ on the number of $\mathbb{Q}$-rational preperiodic points of such a quadratic rational map. In the current work,  we prove the following result.


\begin{Theorem}
\label{thm:intro}
    Let $f:\PP^1\to\PP^1$ be a rational map of degree $2$ defined over $\Q$ with $\Aut(f)\isom\mathfrak{S}_3$. Then $f$ has no rational periodic point of exact period $N>3$.
\end{Theorem} 
In view of Theorem \ref{thm:intro}, we perform a complete analysis of the set of rational preperiodic points $\Preper(f,\Q)$ when $f$ is a rational map of degree $2$ defined over $\Q$ with $\Aut(f) \isom \mathfrak{S}_3$ to obtain the following sharp bound on the number of rational preperiodic points. 

\begin{Theorem}
\label{thm:conj}
Let $f:\PP^1\to\PP^1$ be a rational map of degree $2$ defined over $\Q$ with $\Aut(f) \isom \mathfrak{S}_3$. Then $|\Preper(f,\Q)|\le 6.$
\end{Theorem}

In fact, the bound $6$ in Theorem \ref{thm:conj} is sharp. Theorem \ref{thm:intro} and Theorem \ref{thm:conj} serve as unconditional analogues of Conjecture 2 and Corollary 1 of Poonen, \cite{poonen}, for polynomial maps of degree $2$ defined over $\Q$, and Conjecture 1 and Theorem 2 of Manes, \cite{Manes}, for rational maps of degree $2$ defined over $\Q$ with automorphism group $\mathfrak{C}_2$.

For a map $\phi \in \Q(z)$, the {\em (preperiodic) portrait} $\mathcal P(\phi, \Q)$ is the
directed graph whose vertices are the elements of $\Preper(\phi, \Q)$, with an edge $\alpha\rightarrow \beta$ if and
only if $\phi(\alpha)=\beta$. We classify preperiodic portraits $\mathcal P(f,\mathbb Q)$ when $f$ is a rational map of degree $2$ defined over $\Q$ with $\Aut(f) \isom \mathfrak{S}_3$. In view of the analysis we carry out in \S\ref{sec:different} and \S\ref{sec:preper} on the possible combination of preperiodic points of such maps $f$, we show that the preperiodic portrait $\mathcal P(f,\Q)$ is one of the seven portraits appearing in the following table, see Theorem \ref{thm:classification}. Theorem \ref{thm:intro} is Theorem \ref{thm}, and Theorem \ref{thm:conj} follows from Theorem \ref{thm:classification}. We note that, unlike the case of polynomial maps where $\infty$ is always a fixed point, the set $\Preper(f,\Q)$ may be empty.
\newpage

\renewcommand{\arraystretch}{1} 
\noindent
\begin{tabular}{|p{0.47\textwidth}|p{0.47\textwidth}|}
\hline
\vspace{-0.2cm} \textbf{i.} \par \vspace{0.1cm}
\centering
\begin{tikzcd}[ampersand replacement=\&]
    \bullet \arrow[r] \& \bullet \arrow[out=-45, in=45, loop, distance=2em]
\end{tikzcd}
&
\vspace{-0.2cm} \textbf{ii.} \par \vspace{0.1cm}
\centering
\begin{tikzcd}[ampersand replacement=\&]
    \bullet \arrow[r] \& \bullet \arrow[out=-45, in=45, loop, distance=2em]
\end{tikzcd}
\hspace{0.5cm}
\begin{tikzcd}[ampersand replacement=\&]
    \bullet \arrow[r, bend left] \& \bullet \arrow[l, bend left]
\end{tikzcd}
\tabularnewline[0.6cm] 
\hline

\vspace{-0.2cm} \textbf{iii.} \par \vspace{0.1cm}
\centering
\begin{tikzcd}[ampersand replacement=\&]
    \bullet \arrow[r] \& \bullet \arrow[out=-45, in=45, loop, distance=2em]
\end{tikzcd}
\hspace{0.2cm}
\begin{tikzcd}[ampersand replacement=\&]
    \bullet \arrow[r] \& \bullet \arrow[out=-45, in=45, loop, distance=2em]
\end{tikzcd}
\hspace{0.2cm}
\begin{tikzcd}[ampersand replacement=\&]
    \bullet \arrow[r] \& \bullet \arrow[out=-45, in=45, loop, distance=2em]
\end{tikzcd}
&
\vspace{-0.2cm} \textbf{iv.} \par \vspace{0.1cm}
\centering
\begin{tikzcd}[row sep=small, ampersand replacement=\&]
    \bullet \arrow[r] \& \bullet \arrow[r] \& \bullet \arrow[out=-45, in=45, loop, distance=2em] \\
    \& \bullet \arrow[u] \& 
\end{tikzcd}
\tabularnewline[0.6cm]
\hline

\multicolumn{2}{|p{0.94\textwidth}|}{
    \vspace{-0.2cm} \textbf{v.} \par \vspace{0.1cm}
    \centering
    \begin{tikzcd}[row sep=large, column sep=large, ampersand replacement=\&]
        \& \& \bullet \arrow[dd] \& \bullet \arrow[l] \\
        \bullet \arrow[r] \& \bullet \arrow[ur, bend left=25] \& \& \\
        \& \& \bullet \arrow[ul, bend left=25] \& \bullet \arrow[l]
    \end{tikzcd}
    \vspace*{0.2cm} 
} \tabularnewline[1.2cm] 
\hline

\vspace{-0.2cm} \textbf{vi.} \par \vspace{0.1cm}
\centering
\begin{tikzcd}[ampersand replacement=\&]
    \bullet \arrow[r, bend left] \& \bullet \arrow[l, bend left]
\end{tikzcd}
&
\vspace{-0.2cm} \textbf{vii.} \par \vspace{0.1cm}
\centering
$\emptyset$
\tabularnewline[0.6cm]
\hline

\end{tabular}
\begin{center}
\captionof{table}{Portraits that can be realized as $\mathcal P(f,\mathbb Q)$}
\label{Table1}
\end{center}
 \subsection*{Acknowledgment} The authors are indebted to John R. Doyle for suggesting the proof of Theorem \ref{thm}. This suggestion substantially improved the manuscript, allowing us to replace the conditional results of an earlier draft with the current unconditional results appearing in the current version.
 The authors are grateful to Touko Ursin (Aalto University) for pointing out that the portrait consisting of a single $2$-cycle and the empty portrait were missing from Table \ref{Table1} in an earlier version of this manuscript.
 All the calculations in this work were performed using \Magma, \cite{Magma}, and \Mathematica, \cite{Mathematica}. 
Code verifying the computations made in this paper are available at the GitHub repository \cite{Github}.
 This work builds upon the findings presented in the master's dissertation of H. Bilgili at Sabanc{\i} University under the supervision of M. Sadek. 
  The authors are supported by The Scientific and Technological Research Council of T\"{u}rkiye, T\"{U}B\.{I}TAK, research grant ARDEB 1001/124F352.
  \newpage
\section{Quadratic rational maps with nonabelian automorphism groups} 
Let $K$ be a number field with algebraic closure $\overline K$. 
Given a rational map $f(z)=F(z)/G(z)$, where $F(z),G(z)\in K[z]$ are relatively prime, we set $\deg f=\max\{\deg F,\deg G\}$.
The $n$-th iteration of $f$ is defined by $f^n(z)=f(f^{n-1}(z)),\, n\ge 1$, and $f^0(z)=z.$ Given $a\in \overline K$, the orbit of $a$ under $f$ is the set $\mathcal{O}_f(a)=\{f^n(a),n\ge 0\}.$ A point $a \in \PP^1(K)$ is called {\em periodic} under $f$ if $f^n(a)=a$ for some $n\ge1$. Moreover, if $n$ is the smallest such integer, then $a$ is said to be a periodic point of {\em exact period} $n$; and $\mathcal{O}_f(a)$ is called an {\em $n$-cycle}. 

A point $a\in \PP^1(K)$ is called \textit{preperiodic} under $f$ of type $n_m$, $m\geq 0, n\geq 1$, if $f^{m}(a)$ is periodic of exact period $n$ and, when $m\ge 1$, $f^{m-1}(a)$ is not periodic.
We attach a directed graph to $\mathcal{O}_f(a)$ as follows
$$a\rightarrow f(a)\rightarrow \cdots \rightarrow f^m(a)\rightarrow f^{m+1}(a)\rightarrow\cdots\rightarrow f^{m+n-1}(a)\rightarrow f^{m}(a).$$ 

We set
\[\Preper(f,K)=\{a\in \PP^1(K): a \textrm{ is preperiodic under }f\}.\] 
The set $\Preper(f, K)$ of preperiodic points
of $f$ defined over $K$ can be made a directed graph by drawing an arrow from $a$
to $f(a)$ for each $a\in \Preper(f,K)$.

We define an equivalence relation on rational maps in $K(z)$ of a given degree $d\ge 2$ as follows. Two maps $f_1$ and $f_2$ in $K(z)$ of degree $d\ge2$ are {\em conjugate} if there is $\phi\in \PGL_2(\overline{K})$ such that $f_2=f_1^{\phi}:=\phi\circ f_1\circ \phi^{-1}$. If $\phi\in \PGL_2(K)$, then $f_1$ and $f_2$ are said to be $K$-conjugate. It can be easily seen that if $a$ is a periodic point of exact period $n$ for $f$, then $\phi(a)$ is a point of exact period $n$ for $f^{\phi}$. Moreover, if $f,\phi,$ and $a$ are defined over $K$ such that $f^n(a)=a$, then $g:=f^{\phi}$ and $b:=\phi(a)$ are defined over $K$ with $g^n(b)=b$.
Given $f\in K(z)$, we define
\[\Aut(f)=\{\phi\in\PGL_2(\overline{K}):f^{\phi}=f\}.\]

We will frequently use homogeneous coordinates $f(X,Y)=[F(X,Y),G(X,Y)]$ where $F,G$ are homogeneous polynomials of the same degree with no common factor. It follows that $f^n(X,Y)=[F_n(X,Y),G_n(X,Y)]$ where $F_n$ and $G_n$ are given recursively by 
\begin{eqnarray*} F_n(X,Y)=F_{n-1}(F(X,Y),G(X,Y)),\qquad G_n(X,Y)=G_{n-1}(F(X,Y),G(X,Y)).\end{eqnarray*}
The $n$-period polynomial of $f$ is $\Phi_{f,n}(X,Y)=YF_n(X,Y)-XG_n(X,Y)$. 
\begin{Definition}
The $n$-th {\em dynatomic polynomial of} $f$ is the polynomial
\[\Phi^*_{f,n}(X,Y)=\prod_{e\mid n}\Phi_{f,e}(X,Y)^{\mu(n/e)}\] where $\mu$ is the M\"{o}bius function. 
\end{Definition}
The $n$-th dynatomic polynomial 
$\Phi^*_{f,n}(X,Y)$ is indeed a homogeneous polynomial of degree $\nu_d(n)=\sum_{e\mid n}\mu(n/e)d^e$, where $d=\deg f$, see \cite[Chapter 4, \S 1]{silver}.  
All points of period $e$ under $f$, where $e \mid n$, are roots of the $n$-period polynomial $\Phi_{f,n}$ of $f$. The points of exact period $n$ under $f$ are roots of $\Phi^*_{f,n}$. In fact, the roots of $\Phi^*_{f,n}$ are the points of formal period $n$ under $f$, see  \cite[Chapter 4, \S 1, Theorem 4.5]{silver} for definitions and details. More precisely, the following holds.
\newpage
\begin{Remark}
Let $f (z) \in K(z)$ be a rational map and assume that $P \in \mathbb{P}^1$. 
\begin{enumerate}
    \item If $\Phi_{f,n}(P)=0$, then $P$ is of period $n$ under $f$.
    \item If $\Phi_{f,n}(P)=0$ and $\Phi_{f,m}(P)\neq 0$ for all $m<n$, then  $P$ is a point of exact period $n$ under $f$.
    \item If $P$ is of exact period $n$ under $f$, then $\Phi^*_{f,n}(P)=0$.
\end{enumerate}    
\end{Remark}

The following is \cite[Theorem 5.1]{Milnor}.

\begin{Lemma}
\label{lem:Milnor}
 A quadratic rational map possesses a nontrivial automorphism if and only if it is $\mathbb C$-conjugate to a map in the unique normal form
$f(z) = k(z+1/z)$ with $k\in \mathbb{C}\setminus\{0\}$. For $f$ in this normal form,
$\Aut (f)$ is cyclic of order two (consisting of the maps
$z\mapsto \pm z$) if $k\ne -1/2$, but is nonabelian of order six
for $k =-1/2$.
\end{Lemma}

The classification of preperiodic points of quadratic rational maps with cyclic automorphism groups was discussed in \cite{Manes}. In this work, we focus on quadratic rational maps $f$ over the rational field $\Q$ where $\Aut(f)$ is nonabelian, hence $\Aut(f)\isom \mathfrak{S}_3$. 

If $f$ is a quadratic rational map over a field $K$ with $\Char K\ne 2,3$, then $\Aut(f)\isom\mathfrak{C}_2 $ if and only if $f$ is linearly conjugate over $K$ to a map of the form  $\phi_{k,b}(z)=kz + b/z$ with $k \in K \setminus\{0, -1/2\}$ and $b \in K\setminus\{0\}$. Furthermore, two such maps $\phi_{k,b}$ and $\phi_{k',b'}$ are $K$-conjugate if and only if $k = k'$ and $b/b' \in (K\setminus\{0\})^2$. The map $\phi_{k,b}$ has the automorphism $z\mapsto -z$, see \cite[Lemma 1]{Manes}.

The following is \cite[Lemma 5.1]{Yasuf}.

\begin{Lemma}
\label{lem:Yasuf}
Let $f$ be a quadratic rational map over a field $K$ with $\Char K\ne 2,3$  with $\Aut(f)\isom\mathfrak{S}_3$. Then the map $f$ is $K$-conjugate  to a rational map of the form
\[\theta_{d,k}(z)=\frac{kz^2-2dz+dk}{z^2-2kz+d}\quad\textrm{with}\quad k\in K,\,d\in K\setminus\{0\},\,k^2\ne d.\]
Two maps $\theta_{d,k}$ and $\theta_{d',k'}$ are $K$-conjugate if and only if $$d'=b^2d\quad\textrm{ and } \quad k'\in\left\{\frac{bd}{k},\frac{b(d^2\gamma^3+3dk\gamma^2+3d\gamma+k)}{dk\gamma^3+3d\gamma^2+3k\gamma+1}\right\}
$$ 
for some $\gamma\in  K$ and $b \in K\setminus\{0\} $.
\end{Lemma}

\section{Rational periodic points}
\label{sec:sec}
We now study the arithmetic dynamical behavior of quadratic rational maps over $K$, where $K$ is a number field and most of the time it will be the rational field $\Q$, with automorphism group $\mathfrak{S}_3$. Using Lemma \ref{lem:Yasuf}, any such map is $K$-conjugate to a map of the form $\theta_{d,k}(z)=\frac{kz^2-2dz+dk}{z^2-2kz+d}$ with $k\in K$, $d\in K\setminus \{0\}$, $k^2\ne d$. One also observes that if $k\ne 0$, then conjugating $\theta_{d,k}$ by the scaling $z\mapsto kz$ gives $$\frac{1}{k}\,\theta_{d,k}(kz)=\frac{z^2-2(d/k^2)z+d/k^2}{z^2-2z+d/k^2}=\theta_{d/k^2,1}(z),$$ so $\theta_{d,k}$ is $K$-conjugate to $\theta_{d/k^2,1}$. In the polynomial computations below, this corresponds to the transformation $z\mapsto z/k$, $d\mapsto d/k^2$.  
\newpage
\begin{Lemma}
\label{lem:sameimage}
Let $K$ be a field of characteristic different from $2$ and $3$. Let $\theta_{d,k}(z)=(kz^2-2dz+dk)/(z^2-2kz+d)$ with $k\in K$, $d\in K\setminus \{0\}$, $k^2\ne d$. 
\begin{itemize}
\item[i)] Let $z_1,z_2\in \PP^1(K)$ be distinct. One has $\theta_{d,k}(z_1)=\theta_{d,k}(z_2)$ if and only if $z_2=d/z_1$, where we use the convention $d/0=\infty$ and $d/\infty=0$.
\item[ii)] If $z$ is a $K$-rational periodic point of $\theta_{d,k}$ such that $d\ne z^2$, then $d/z\in\Preper(\theta_{d,k},K)$. In particular, 
if $z$ is a $K$-rational fixed point of $\theta_{d,k}$ such that $d\ne z^2$, then $d/z$ is a $K$-rational preperiodic point of type $1_1$. 
\end{itemize}
\end{Lemma}
\begin{Proof}
  i) Since the resultant of $kz^2-2dz+dk$ and $z^2-2kz+d$ is $4d(d-k^2)^2\ne 0$, one has $\theta_{d,k}(z)=\infty$ for $z\in K$ if and only if $z^2-2kz+d=0$. The two roots of the latter polynomial are distinct as $k^2\ne d$, and their product is $d$; this proves the statement when $\theta_{d,k}(z_1)=\infty$. For $z_1,z_2\in K$ with $\theta_{d,k}(z_1),\theta_{d,k}(z_2)\ne\infty$, one has
  \[\theta_{d,k}(z_1)-\theta_{d,k}(z_2)=\frac{2(d-k^2)(z_1-z_2)(z_1z_2-d)}{(z_1^2-2kz_1+d)(z_2^2-2kz_2+d)},\]
  and the statement follows since $d\ne k^2$. Finally, $\theta_{d,k}(0)=k=\theta_{d,k}(\infty)$, and for $z\in K$ one has $\theta_{d,k}(z)-k=2z(k^2-d)/(z^2-2kz+d)$, so $\theta_{d,k}(z)=k$ if and only if $z=0$. ii) follows immediately from i).
\end{Proof}

\begin{Proposition}
\label{prop:period1}
Let $K$ be a field of characteristic different from $2$ and $3$. Let $\theta_{d,k}(z)=(kz^2-2dz+dk)/(z^2-2kz+d)$ with $k\in K$, $d\in K\setminus \{0\}$, $k^2\ne d$. 

If $k=0$, then $0$ is always a fixed point of $\theta_{d,0}$. If $k\ne 0$, then $\theta_{d,k}$ has a $K$-rational fixed point if and only if $$k=u,\qquad d=-\frac{(3 u - v) v^2}{ u - 3 v},\quad u,v\in K,\, u\ne 3v,$$ where the fixed point is $z=v$. 
\end{Proposition}
\begin{Proof}
One sees that 
 \begin{equation}
     \begin{split}
         \Phi_ {\theta_{d,k},1} (z)=
    d k - 3 d z + 3 k z^2 -  z^3. \label{eq1}
     \end{split}
 \end{equation}
 It is clear that when $k=0$, one has that $z=0$ is a fixed point of $\theta_{d,0}$. So, we now assume that $k\ne 0$. Assuming that $z=v$ is a $K$-rational fixed point for $\theta_{d,k}$, the result now follows by solving for $d$ in terms of $z=v$ and $k=u$. 
\end{Proof}

\begin{Proposition}
\label{prop:period2}
Let $K$ be a field of characteristic different from $2$ and $3$. Let $\theta_{d,k}(z)=(kz^2-2dz+dk)/(z^2-2kz+d)$ with $k\in K$, $d\in K\setminus \{0\}$, $k^2\ne d$. Then $\theta_{d,k}$ has a periodic point of exact period $2$ if and only if $d=c^2$, $c\in K$. In this case, there are exactly two such points $c,-c$, and they form a $2$-cycle. 
\end{Proposition}
\begin{Proof}
   This follows directly from the following equality
   \begin{equation*}
     \begin{split}
         \Phi^{*}_ {\theta_{d,k},2} (z) &=  
           -\left(\left(d-k^2\right) \left(d-z^2\right)\right).
     \end{split}
 \end{equation*} 
\end{Proof}
\begin{Remark}
\label{rem:2cyclenotail}
Assume that $d=c^2$, $c\in K$. Since $\theta_{c^2,k}(\mp c)=\pm c$ and $d/(\mp c)=\mp c$, Lemma \ref{lem:sameimage} i) implies that $\mp c$ is the only point $z\in\PP^1(K)$ with $\theta_{c^2,k}(z)=\pm c$. In particular, $\theta_{d,k}$ possesses no $K$-rational preperiodic point of type $2_m$, $m\ge 1$. Note also that $\theta_{c^2,k}(\infty)=k\ne\pm c$ since $k^2\ne d$.
\end{Remark}

\begin{Theorem}
\label{thm:period3}
 Let $\theta_{d,k}(z)=(kz^2-2dz+dk)/(z^2-2kz+d)$ with $k\in \Q$, $d\in \Q\setminus \{0\}$, $k^2\ne d$. 
 
 The map $\theta_{d,0}$ has no rational periodic points with exact period $3$; whereas the map $\theta_{d,k}$, $k\ne0$, has a periodic point $x\in\Q$ of exact period $3$ if and only if $d=-u^2/3$ and $k$ is one of the two values $\displaystyle\frac{\pm u^4 + 3 u^3 x\mp 9 u^2 x^2 - 3 u x^3}{
u^3\mp 9 u^2 x - 9 u x^2\pm 9 x^3}$, for some $u\in \Q\setminus\{0\}$. In particular, $-3d$ is a rational square. In this case, the $3$-cycle is given by $$x\longrightarrow\frac{u (\pm u + x)}{u \mp 3 x}\longrightarrow \mp\frac{u (u \mp x)}{u \pm 3 x}\longrightarrow x,$$ where a point of the $3$-cycle whose denominator vanishes, i.e., when $u=\pm 3x$, is understood to be $\infty$.
In addition, for every point $z$ in the $3$-cycle, $d/z$ is a $\Q$-rational preperiodic point of type $3_1$.
\end{Theorem}
\begin{Proof}
If $k=0$, then a map $\theta_{d,0}$ together with a rational periodic point of exact period $3$ will give rise to a rational root of the polynomial $3 d^3 + 27 d^2 z^2 + 33 d z^4 + z^6$. Since the polynomial $3 + 27 Z + 33  Z^2 + Z^3$, $Z=z^2/d$, is irreducible over $\Q$, we conclude that there is no $d\in \Q$ such that $\theta_{d,0}$ possesses a rational periodic point with exact period $3$.

 We now assume $k\ne0$. We compute the $3$-rd dynatomic polynomial of $\theta_{d,k}$ to obtain
     {\footnotesize\begin{equation*}
   \begin{split}
       \Phi^{*}_ {\theta_{d,k},3} (z)
        &= (d-k^2)^2 (3 d^4 + 27 d^3 z^2 + 33 d^2 z^4 + d z^6 + 
 k (-24 d^3 z - 80 d^2 z^3 - 24 d z^5) + 
 k^2 (d^3 + 33 d^2 z^2 + 27 d z^4 + 3 z^6)).
    \end{split}
    \end{equation*}}
         Recall that if $\theta_{d,k}$ possesses a rational periodic point $x$ of exact period $3$, then $ \Phi^{*}_ {\theta_{d,k},3} (x)=0$. One sees that $ \Phi^{*}_ {\theta_{d,k},3}(z)/(d-k^2)^2 $ is a polynomial of degree $2$ in $k$; its leading coefficient $d^3+33d^2z^2+27dz^4+3z^6$ does not vanish for $z\in\Q$, since $3Z^3+27Z^2+33Z+1$, $Z=z^2/d$, has no rational root. Its roots are
     \begin{align*}
          \frac{12 d^3 z+40 d^2 z^3+12 d z^5\pm \sqrt{-3d \left(d-z^2\right)^6}}{d^3+33 d^2 z^2+27 d z^4+3 z^6}.
         \end{align*}
        If $d\ne z^2$, the latter roots are rational if and only if $-3d$ is a rational square. If $d=z^2$, then the two roots coincide and are equal to $z$, so $k^2=d$, which is excluded. Computing the points in the $3$-cycle is straightforward. Finally, the fact that $-3d$ must be a rational square together with Lemma \ref{lem:sameimage} imply that for every point $z$ in the $3$-cycle, $d/z$ is a $\Q$-rational preperiodic point of type $3_1$. Computations can be found in \href{https://github.com/hasanbilgili44/Calculation-for-Paper/blob/main/Nonabelian-Automorphism/Case%20Period%203}{Case Period 3} of \cite{Github}.
\end{Proof}
We remark that the denominators $u^3\mp 9u^2x-9ux^2\pm 9x^3$ in Theorem \ref{thm:period3} do not vanish for $u\in\Q\setminus\{0\}$, since the polynomials $t^3\mp 9t^2-9t\pm 9$ have no rational roots.

For $\theta_{d,k}$, $k\ne0$, one sees that $\theta_{d,k}(\infty)=k$ and $\theta_{d,k}(k)=-k$. Hence $\infty$ is a periodic point of exact period $3$ if and only if $\theta_{d,k}(-k)=\infty$, i.e., $d=-3k^2$. In this case, the $3$-cycle is $\infty\rightarrow k\rightarrow -k\rightarrow \infty$; in Theorem \ref{thm:period3}, this corresponds to $u=\pm 3x$.

In Theorem \ref{thm:period3}, if $k=x$, then $\infty$ is a preperiodic point. For $k=\displaystyle\frac{ u^4 + 3 u^3 x- 9 u^2 x^2 - 3 u x^3}{
u^3- 9 u^2 x - 9 u x^2+ 9 x^3}$, the equality $k=x$ implies that $u=x$ or $u=-3x$. If $u=-3x$, then $\infty$ lies in the $3$-cycle $\infty\rightarrow x\rightarrow -x\rightarrow\infty$. If $u=x$, then $\infty$ is a preperiodic point of type $3_1$ with the orbit $\infty\rightarrow u\rightarrow -u\rightarrow 0\rightarrow u$. Similarly, if $k=\displaystyle\frac{- u^4 + 3 u^3 x+ 9 u^2 x^2 - 3 u x^3}{
u^3+ 9 u^2 x - 9 u x^2- 9 x^3}$, then $k=x$ implies that $u=-x$ or $u=3x$. If $u=3x$, then $\infty$ lies in the $3$-cycle $\infty\rightarrow x\rightarrow -x\rightarrow\infty$. If $u=-x$, then $\theta_{d,k}$ possesses the orbit $\infty\rightarrow -u\rightarrow u\rightarrow 0\rightarrow -u$.

\section{Non-existence of rational periodic points of periods $\ge 4$}
In this section, we prove that the map $\theta_{d,k}(z)=(kz^2-2dz+dk)/(z^2-2kz+d)$ with $k\in \Q$, $d\in \Q\setminus \{0\}$, $k^2\ne d$, cannot possess a rational periodic point of exact period exceeding $3$.

For a quadratic rational map $f$ over $K$, 
there is one $\overline{K}$-conjugacy class of maps with automorphism group $\mathfrak{S}_3$. In particular, $1/z^2$ lies in the latter conjugacy class with automorphism group generated by $z \mapsto 1/z$ and $z\mapsto \omega z$
where $\omega$ is a primitive cube root of unity. 

\begin{Theorem}
\label{thm}
Let $f:\PP^1\to\PP^1$ be a rational map of degree $2$ defined over $\Q$ with $\Aut(f)\isom\mathfrak{S}_3$. Then $f$ has no rational periodic point of exact period $N>3$.   
\end{Theorem}
\begin{Proof}
To obtain a contradiction, assume that $f$ admits a cycle of rational points of length $n \ge 4$. The only degree-2 rational map with nonabelian automorphism group is $g(z) := 1/z^2$, up to conjugation over $\overline{\mathbb{Q}}$, see \cite{Yasuf}. Thus, there exists an element of $\PGL_2(\overline{\mathbb{Q}})$, say $\gamma$, such that $f = \gamma \circ g \circ \gamma^{-1}$. Moreover, $\gamma$ must map any of the $n$-cycles for $g$ to an $n$-cycle for $f$.

By the assumption $n \ge 4$, one can easily see that the $n$-cycle for $g$ consists of roots of unity; in fact, an orbit of $g$ will have the following form 
\[
\zeta \mapsto \zeta^{-2} \mapsto \zeta^4 \mapsto \zeta^{-8} \mapsto \cdots \mapsto \zeta^{(-2)^n} = \zeta
\]
for some (odd-order) root of unity $\zeta$.

Since $\gamma(\zeta)$, $\gamma(\zeta^{-2})$, and $\gamma(\zeta^4)$ are rational, we can find an element $\gamma'$ of $\PGL_2(\mathbb{Q})$ that maps $\gamma(\zeta)$, $\gamma(\zeta^{-2})$, and $\gamma(\zeta^4)$ to $0$, $\infty$, and $1$, respectively. Define $\eta$ as $\gamma' \circ \gamma$. Then, we have
\[
\eta(\zeta) = 0, \quad \eta(\zeta^{-2}) = \infty, \quad \text{and} \quad \eta(\zeta^4) = 1.
\]
This uniquely determines $\eta$; specifically, we have
\[
\eta(z) = \frac{z-\zeta}{z-\zeta^{-2}} \cdot \frac{\zeta^4-\zeta^{-2}}{\zeta^4-\zeta}.
\]
By our assumption, the $n$-cycle for $f$ consists of rational points. This leads us to the fact that $\gamma(\zeta^{-8}) \in \mathbb{Q}$. So, $\eta(\zeta^{-8}) \in \mathbb{Q}$. Then, we have
\begin{align*}
\eta(\zeta^{-8}) = \frac{\zeta^{-8}-\zeta}{\zeta^{-8}-\zeta^{-2}} \cdot \frac{\zeta^4-\zeta^{-2}}{\zeta^4-\zeta} 
= \zeta^3 + 1 + \zeta^{-3}.
\end{align*}
Recall that $\zeta$ has odd order. Thus, $\zeta^2$ is a $\mathrm{Gal}(\overline{\mathbb{Q}}/\mathbb{Q})$-conjugate of $\zeta$. Since $\eta(\zeta^{-8})$ is rational, we have
\[
\zeta^6 + \zeta^{-6} = \zeta^3 + \zeta^{-3}.
\]
Nevertheless, this equality holds only if $\zeta$ is a 9th root of unity. Under the map $1/z^2$, such a root has a maximum period of 3, which directly contradicts the assumption that $n \ge 4$.
\end{Proof}

\section{Rational maps with different cycles}
\label{sec:different}
In this section, we discuss the possibilities for the map $\theta_{d,k}$ to have at least two distinct cycles. Equivalently, we investigate whether there are rational periodic points $z_1,\cdots,z_m$ under $\theta_{d,k}$, $m\ge 2$, such that $\mathcal{O}_{\theta_{d,k}}(z_1),\cdots,\mathcal{O}_{\theta_{d,k}}(z_m)$ are disjoint.  
\begin{Proposition}
    \label{prop:period1multiples}
Let $K$ be a field of characteristic different from $2$ and $3$. Let $\theta_{d,k}(z)=(kz^2-2dz+dk)/(z^2-2kz+d)$ with $k\in K$, $d\in K\setminus \{0\}$, $k^2\ne d$. 

 If $k=0$, then either $\theta_{d,0}$ possesses three fixed points, $0$ and $\pm \sqrt{-3d}$, when $-3d$ is a $K$-rational square; or  $0$ is the only fixed point of $\theta_{d,0}$ in $K$, when $-3d$ is not a $K$-rational square. 

If $k\ne0$, then $\theta_{d,k}$ has three fixed points if and only if $$ k=\frac{-2 (3 s - 5 t) t}{3 (s - 3 t) (3 s - t)},\quad d=\frac{-(3 s - 13 t)^2 (3 s - 7 t)^2}{2187 (s - 3 t)^4}, \quad s,t\in K, \,s\not\in\{3t,t/3\},$$
 where the three fixed points are 
 \[\frac{(3 s - 13 t) (3 s - 7 t)}{9 (s - 3 t) (3 s - t)},\quad \frac{-2 (3 s - 13 t) t}{27 (s - 3 t)^2},\quad \frac{-(3 s - 7 t) (3 s - 5 t)}{27 (s - 3 t)^2}.\]
\end{Proposition}
\begin{Proof}
   When $k=0$, this is straightforward by computing $\Phi_{\theta_{d,0},1}$. 

   We assume $k\ne0$. Using Proposition \ref{prop:period1}, we know that the existence of one fixed point imposes a parametrization of $d$ and $k$ in terms of two $K$-rational parameters $u$ and $v$. Computing $\Phi_{\theta_{d,k},1}$ implies that the other two fixed points of $\theta_{d,k}$ are $K$-rational if and only if $3(u - 3 v)  (3 u - v)$ is a $K$-rational square. In particular, one has
   \[u = \frac{-2 (3 s - 5 t) t}{3 (s - 3 t) (3 s - t)},\qquad 
v = \frac{(3 s - 13 t) (3 s - 7 t)}{9 (s - 3 t) (3 s - t)},\quad s\not\in\{3t,t/3\},\] hence the result follows.
\end{Proof}

\begin{Remark}
 \label{rem:morecycles0} In view of Proposition \ref{prop:period2},  Theorem \ref{thm:period3}, and Proposition \ref{prop:period1multiples}, we recall that the map $\theta_{d,0}(z)=(-2dz)/(z^2+d)$, $d\in \Q\setminus \{0\}$ satisfies the following statements.
 \begin{itemize}
     \item [i)] $\theta_{d,0}$ possesses exactly one $\Q$-rational fixed point, namely $0$, and a $\Q$-rational periodic point of exact period $2$ if and only if $d $ is a $\Q$-rational square.
     \item[ii)] If $\theta_{d,0}$ possesses three $\Q$-rational fixed points, then it contains no $\Q$-rational periodic points of exact period $2$.  
     \item[iii)] $\theta_{d,0}$ has no $\Q$-rational periodic points of exact period $3$. 
 \end{itemize}
\end{Remark}

\begin{Theorem}
\label{thm:morecycles}
 Let $\theta_{d,k}(z)=(kz^2-2dz+dk)/(z^2-2kz+d)$ with $k\in \Q$, $d\in \Q\setminus \{0\}$, $k^2\ne d$. 
\begin{itemize}
\item[i)] The map $\theta_{d,k}$ possesses exactly one $\Q$-rational fixed point and a $\Q$-rational periodic point of exact period $2$ if and only if $$k=\frac{3 s^2 - 6 s t + 7 t^2}{2 (s - 3 t) (3 s - t)},\qquad d=\frac{(9 s^2 - 42 s t + 61 t^2)^2}{324 (s - 3 t)^4}, \quad s,t\in\Q, s\not\in\{3t,t/3\},$$ with the fixed point $(9 s^2 - 42 s t + 61 t^2)/(6 (s - 3 t) (3 s - t))$ and the  $2$-cycle consisting of $\sqrt{d}$ and $-\sqrt{d}$. 
\item[ii)] The map $\theta_{d,k}$ cannot possess a $\Q$-rational fixed point and a $\Q$-rational periodic point of exact period $3$.
\item[iii)] If $\theta_{d,k}$ possesses three $\Q$-rational fixed points, then it has no $\Q$-rational periodic points of either exact periods $2$ or $3$.
\item[iv)] The map $\theta_{d,k}$ cannot have more than $m$ $\Q$-rational periodic points of exact period $m$, $m=2,3$.
\item[v)]  The map $\theta_{d,k}$ cannot have a $\Q$-rational periodic point of exact period $2$ together with a $\Q$-rational periodic point of exact period $3$.
\end{itemize}
\end{Theorem}
\begin{Proof}
i) In view of Proposition \ref{prop:period1} and Proposition \ref{prop:period2}, such a map is of the form $\theta_{d,k}$ with $k=u$ and $ d=-(3u-v)v^2/(u-3v)$, moreover, it must satisfy that $-(u-3v)(3u-v)$ is a $\Q$-rational square. It follows that 
\begin{equation*}
    u=\frac{3 s^2 - 6 s t + 7 t^2}{2 (s - 3 t) (3 s - t)},\quad 
v=\frac{9 s^2 - 42 s t + 61 t^2}{6 (s - 3 t) (3 s - t)}.
\end{equation*}

Substituting the latter expressions in $k$ and $d$ yields the result. 

ii) If $k=0$, then $\theta_{d,0}$ has no $\Q$-rational periodic point of exact period $3$, see Theorem \ref{thm:period3}. Assume $k\ne 0$. According to Proposition \ref{prop:period1}, the map $\theta_{d,k}$ has a $\Q$-rational fixed point if $k=u$ and $d$ is of the form $-(3 u - v) v^2/( u - 3 v)$. Then, one obtains that
{\footnotesize\begin{equation*}    
\begin{split}
 \Phi^*_{\theta_{d,k},3}(z) &= \frac{1}{(u-3 v)^4} \Big(-3 u^6 \left(9 v^6-99 v^4 z^2+27 v^2 z^4-z^6\right)\\ &+36 u^5 v \left(3 v^6+18 v^5 z-55 v^4 z^2-20 v^3 z^3+21 v^2 z^4+2 v z^5-z^6\right)  \\ 
 &+  3 u^4 v^2 \left(51 v^6-864 v^5 z+1055 v^4 z^2+1600 v^3 z^3-711 v^2 z^4-224 v z^5+53 z^6\right) \\ &-8 u^3 v^3 \left(37 v^6-270 v^5 z-117 v^4 z^2+1180 v^3 z^3-117 v^2 z^4-270 v z^5+37 z^6\right) \\ &+3 u^2 v^4 \left(53 v^6-224 v^5 z-711 v^4 z^2+1600 v^3 z^3+1055 v^2 z^4-864 v z^5+51 z^6\right) \\ &-36 u v^5 \left(v^6-2 v^5 z-21 v^4 z^2+20 v^3 z^3+55 v^2 z^4-18 v z^5-3 z^6\right) \\ &+3 v^6 \left(v^6-27 v^4 z^2+99 v^2 z^4-9 z^6\right) \Big).
\end{split}  
\end{equation*}}
The vanishing of the homogeneous polynomial $(u-3 v)^4\Phi^*_{\theta_{d,k},3}(z)$ defines a planar curve $C$ that is reducible 
   over the number field defined by $f(x) :=x^3 + 34992x^2 + 76527504x -
    2752235153856$, \Magma, \cite{Magma}. One may also see that the curve $C$ consists of the three irreducible components
    {\footnotesize\allowdisplaybreaks
    \begin{align*}
            C_1 &: z^2u^2 + 1/337996476(-\alpha^2 - 56862\alpha - 1658095920)z^2uv +1/112665492(\alpha^2 + 56862\alpha + 306110016)zu^2v 
            \\ &+ 1/112665492(\alpha^2+ 56862\alpha + 644106492)z^2v^2 
            + 1/168998238(-5\alpha^2 - 284310\alpha -1530550080)zuv^2 
            \\ &+ 1/112665492(\alpha^2 + 56862\alpha +644106492)u^2v^2 
            + 1/112665492(\alpha^2 + 56862\alpha + 306110016)zv^3 
            \\ &+1/337996476(-\alpha^2 - 56862\alpha - 1658095920)uv^3 + v^4=0, \\
            C_2&:z^2u^2 + 1/675992952(-5\alpha^2 - 52488\alpha - 1530550080)z^2uv
            +1/225330984(5\alpha^2 + 52488\alpha - 1173421728)zu^2v 
            \\ &+1/225330984(5\alpha^2 + 52488\alpha - 497428776)z^2v^2 
            +1/337996476(-25\alpha^2 - 262440\alpha + 5867108640)zuv^2 
            \\ &+1/225330984(5\alpha^2 + 52488\alpha - 497428776)u^2v^2 
            +1/225330984(5\alpha^2 + 52488\alpha - 1173421728)zv^3 
            \\ &+1/675992952(-5\alpha^2 - 52488\alpha - 1530550080)uv^3 + v^4=0, \\
            C_3&: z^2u^2 + 1/675992952(7\alpha^2 + 166212\alpha - 3265173504)z^2uv 
            +1/225330984(-7\alpha^2 - 166212\alpha + 561201696)zu^2v 
            \\ &+1/225330984(-7\alpha^2 - 166212\alpha + 1237194648)z^2v^2 
            +1/337996476(35\alpha^2 + 831060\alpha - 2806008480)zuv^2 
            \\ &+1/225330984(-7\alpha^2 - 166212\alpha + 1237194648)u^2v^2 
            +1/225330984(-7\alpha^2 - 166212\alpha + 561201696)zv^3 
            \\ &+1/675992952(7\alpha^2 + 166212\alpha - 3265173504)uv^3 + v^4=0,
    \end{align*}}%
    where $\alpha$ is a root of $f(x)$. The curve $C_i$ can be described by $g_i(z,u,v)+\alpha h_i(z,u,v)+\alpha^2l_i(z,u,v)=0$ where $g_i, h_i$ and $l_i$ are elements of $\mathbb{Q}[z,u,v]$, for $i=1,2,3$. The linear independence of $1$, $\alpha$ and $\alpha^2$ over $\mathbb{Q}$ implies that a $\mathbb{Q}$-rational point on the curve defined by $(u-3 v)^4\Phi^*_{\theta_{d,k},3}(z)=0$ provides a rational solution to at least one of the following system of equations 
    \begin{equation*}
        g_i(z,u,v)= h_i(z,u,v)=l_i(z,u,v)=0, \ \ i=1,2,3.
    \end{equation*}
    Magma asserts that the intersection of the planar curves given by $g_i(z,u,v)= h_i(z,u,v)=l_i(z,u,v)=0$, $i=1,2,3$, consists of the following points
    \begin{equation*}
        (z,u,v)\in \{ (0 : 1/3 : 1), (1 : 1 : 1), (0 : 1 : 0), (1 : 0 : 0)  \}.
    \end{equation*}
    In fact, the latter points are the singular points of the curve $C$ and none of them gives rise to a quadratic
rational map $\theta_{d,k}$, $d\ne 0$, $k^2\ne d$. The detailed computations can be found in \href{https://github.com/hasanbilgili44/Calculation-for-Paper/blob/main/Nonabelian-Automorphism/Theorem%205.3}{Theorem 5.3} of \cite{Github}. 

iii) According to Proposition \ref{prop:period1multiples}, the existence of three $\Q$-rational fixed points implies that $-3d$ is a $\Q$-rational square which contradicts the fact that $d$ must be a $\Q$-rational square when $\theta_{d,k}$ possesses a $\Q$-rational periodic point of exact period $2$, see Proposition \ref{prop:period2}. That $\theta_{d,k}$ cannot possess three $\Q$-rational fixed points together with a $\Q$-rational periodic point of exact period $3$ follows from ii).

iv) follows right away for $m=2$ from the fact that the $2$-nd dynatomic polynomial is of degree $2$. As for $m=3$, Theorem \ref{thm:period3} yields that $x,u(\pm u+x)/(u\mp 3x), \mp u(u\mp x)/(u\pm 3x)$ form a $3$-cycle for $\theta_{d,k}$. Replacing $u$ by $-u$ interchanges the two signs in Theorem \ref{thm:period3}, so it suffices to consider the lower signs. Therefore, if $\theta_{d,k}$ possesses another $3$-cycle, then this will give rise to a rational point on the curve 
{\footnotesize\begin{eqnarray*}
C:\,u^6 + 3 u^5 x - 9 u^4 x^2 - 3 u^3 x^3 - 3 u^5 z + 27 u^4 x z+ 
  27 u^3 x^2 z - 27 u^2 x^3 z&-& 9 u^4 z^2 - 27 u^3 x z^2 + 
  81 u^2 x^2 z^2 + 27 u x^3 z^2\\&+& 3 u^3 z^3 - 27 u^2 x z^3 - 
  27 u x^2 z^3 + 27 x^3 z^3=0.
  \end{eqnarray*}}
  The curve $C$ has three irreducible components over the number field $\Q(\alpha)$ where $\alpha$ is a root of the polynomial $f(t)=t^3 - 2187t^2 + 1358127t -
    176969853$. Using the linear independence of $1,\alpha,\alpha^2$, one verifies that the only $\Q$-rational points on any of these irreducible components are the singular points $(z:x:u)=(0:1:0)$, $(1:0:0)$, hence the result for $m=3$. For the details of the computations, see \href{https://github.com/hasanbilgili44/Calculation-for-Paper/blob/main/Nonabelian-Automorphism/Theorem%205.3}{Theorem 5.3} of \cite{Github}. This argument also covers the case $u=\pm 3x$, in which the $3$-cycle contains $\infty$. For instance, if $d=-3k^2$, then $\theta_{d,k}$ is conjugate to $\theta_{-3,1}$, and up to a nonzero constant factor one has $\Phi^{*}_{\theta_{-3,1},3}(z)=(z-1)(z+1)(z^3+3z^2-9z-3)$, where $\infty$ is the remaining root of the homogenization of degree $6$. Since $z^3+3z^2-9z-3$ has no rational roots, $\theta_{-3,1}$ possesses exactly one $\Q$-rational $3$-cycle, namely $\infty\rightarrow 1\rightarrow -1\rightarrow\infty$.

    v) follows in a similar manner from Proposition \ref{prop:period2} and Theorem \ref{thm:period3}.
\end{Proof}
\section{Rational preperiodic points}
\label{sec:preper}
In this section, we study the $\Q$-rational preperiodic points of $\theta_{d,k}$.
For $m, n\ge  1$, we define the generalized $(m, n)$-dynatomic
polynomial of a rational map $\phi$ to be
\begin{equation*}
    \Phi^{*}_{\phi,m,n}(X,Y)=\frac{\Phi^{*}_{\phi,n}(F_{m}(X,Y),G_{m}(X,Y))}{\Phi^{*}_{\phi,n}(F_{m-1}(X,Y),G_{m-1}(X,Y))}.
\end{equation*}
If $P$ is a preperiodic point under $\phi$ of type $n_m$, then $\Phi^{*}_{\phi,m,n}(P)=0$, see \cite[Chapter 4, Exercise 4.11]{silver}.
\begin{Theorem}
\label{thm:preper0}
     Consider the map $\theta_{d,0}(z)=-2dz/(z^2+d)$, where $d\in \Q\setminus \{0\}$.
    \begin{itemize}
        \item [i)] $\theta_{d,0}$ possesses a $\Q$-rational preperiodic point (different from $\infty$) of type $1_1$ if and only if $-3d$ is a $\Q$-rational square. In this case, $\theta_{d,0}$ possesses the following three orbits 
        \begin{eqnarray*}
         \pm a \rightarrow \mp 3a \rightarrow \mp 3a,\qquad \infty\rightarrow 0\rightarrow 0,
     \end{eqnarray*}
        where $d=-3a^2$, $a\in \Q\setminus\{0\}$. 
        \item[ii)] $\theta_{d,0}$ possesses a $\Q$-rational preperiodic point of type $1_2$ if and only if $-d$ is a $\Q$-rational square. In this case, $\theta_{d,0}$ possesses the following orbit 
        \begin{eqnarray*}
         \pm a \rightarrow \infty \rightarrow 0 \rightarrow 0,
     \end{eqnarray*}
     where $d=-a^2$, $a\in \Q\setminus\{0\}$
        \item[iii)]  $\theta_{d,0}$ possesses no $\Q$-rational preperiodic point of type $1_m$, $m\ge 3$ 
    \end{itemize}
\end{Theorem}
\begin{Proof}
     i) This follows immediately from the fact that $\Phi^{*}_{\theta_{d,0},1,1}(z)= -2 d (d^2 + 3 d z^2)$ whose vanishing is equivalent to  $-3d$ being a $\Q$-rational square. One may check that $\theta_{d,0}(\pm a)=\mp3a$, and $\theta_{d,0}(\mp 3a)=\mp 3a$.
   
     ii) The generalized $(2,1)$-dynatomic polynomial of $\theta_{d,0}$ is given by
    \begin{equation*}
        \Phi^{*}_{\theta_{d,0},2,1}(z) = -2 d^5 (z^2 +d) \left(\frac{z^4}{d^2} + 14 \frac{z^2}{d} +1\right).
    \end{equation*}
    Since $Z^2+14Z+1$ is irreducible over $\Q$, and $d$ cannot be $0$, the vanishing of $\Phi^{*}_{\theta_{d,0},2,1}(z)$ is equivalent to $-d$ being a $\Q$-rational square. 
     
     iii) The generalized $(3,1)$-dynatomic polynomial of $\theta_{d,0}$ is defined by
     \begin{equation*}
     \begin{split}
          \Phi^{*}_{\theta_{d,0},3,1}(z)&=  -2 d^4 (d^8 + 66 d^7 z^2 + 495 d^6 z^4 + 924 d^5 z^6 + 495 d^4 z^8 + 
   66 d^3 z^{10} + d^2 z^{12}) \\
   &=-2 d^6 (z^4 + 6z^2d + d^2) (z^8 + 60z^6d + 134z^4d^2 + 60z^2d^3 + d^4).
     \end{split}
     \end{equation*}
     Applying the transformation $z \mapsto \frac{z^2}{d}$, the polynomial $\Phi^{*}_{\theta_{d,0},3,1}(z)$ possesses two irreducible factors over $\mathbb{Q}$, hence, iii) holds. The relevant computations can be found in \href{https://github.com/hasanbilgili44/Calculation-for-Paper/blob/main/Nonabelian-Automorphism/Theorem%206.1}{Theorem 6.1} of \cite{Github}. 
    
\end{Proof}
\begin{Theorem}
\label{thm:preper1}
Let $\theta_{d,k}(z)=(kz^2-2dz+dk)/(z^2-2kz+d)$ with $k\in \Q$, $d\in \Q\setminus \{0\}$, $k^2\ne d$. 
\begin{itemize}
\item[i)] If $z$ is a $\Q$-rational fixed point of $\theta_{d,k}$, then $d\ne z^2$ and $d/z$ is a $\Q$-rational preperiodic point of type $1_1$. 
    \item[ii)] $\theta_{d,k}$, $k\ne 0$, possesses a $\Q$-rational preperiodic point $z$ of type $1_2$ if and only if {\footnotesize$$d=-\frac{k^2(8 s - 41 t)^2 (64 s^2 - 272 s t + 97 t^2)^2}{(8 s - 
    25 t)^2 (64 s^2 - 784 s t + 2209 t^2)^2},\quad   z=\frac{k(8 s - 41 t) (8 s - 33 t) (64 s^2 - 272 s t + 97 t^2)}{8 (8 s - 
   25 t) t (64 s^2 - 784 s t + 2209 t^2)}, \quad s,t\in\Q, t\ne0,s\ne 25t/8.$$} In this case, $\theta_{d,k}$ has only one $\Q$-rational fixed point.
    \item[iii)] $\theta_{d,k}$ possesses no $\Q$-rational preperiodic point of type $1_m$, $m\ge 3$.
\end{itemize}
\end{Theorem}
\begin{Proof}
    i) 
  In accordance with Lemma \ref{lem:sameimage}, a $\Q$-rational fixed point $z$ of $\theta_{d,k}$ is the image of a $\Q$-rational preperiodic point of type $1_1$ if and only if $d\ne z^2$. For $\theta_{d,0}$, the statement follows as $d\ne 0$ and $-3d$ is a $\Q$-rational square when $\theta_{d,0}$ has three $\Q$-rational fixed points, see Proposition \ref{prop:period1multiples}. If $\theta_{d,k}$, $k\ne0$, has exactly one $\Q$-rational fixed point $z$, then setting $d=z^2$ in Proposition \ref{prop:period1} yields that $u$ must be the same as $v$. However, this implies that $\theta_{d,k}$ is a constant map. In Proposition \ref{prop:period1multiples}, one may easily check that for any $s,t\in\Q$, $d\ne z^2$ for any of the three $\Q$-rational fixed points $z$ of $\theta_{d,k}$.

ii) The generalized $(2,1)$-dynatomic polynomial of $\theta_{d,k}$ is defined by 
\begin{equation*}
         \Phi^{*}_{\theta_{d,k},2,1}(z)=-2(d-k^2)^2\big(d^4 + k^2 z^6 + d z^4 (15 k^2 - 12 k z + z^2) + 
    5 d^2 z^2 (3 k^2 - 8 k z + 3 z^2) + d^3 (k^2 - 12 k z + 15 z^2)\big).
\end{equation*} 
Applying the transformation $z\mapsto z/k, d\mapsto d/k^2$, the vanishing of $\Phi^{*}_{\theta_{d,k},2,1}(z)$ yields a rational point on the affine curve defined by
\[d^3 + d^4 - 12 d^3 z + 15 d^2 z^2 + 15 d^3 z^2 - 40 d^2 z^3 + 
  15 d z^4 + 15 d^2 z^4 - 12 d z^5 + z^6 + d z^6=0.\]
  The latter is a genus $0$ curve with three singular points $(d:z:w)=(0 : 0 : 1), (1 : 1 : 1), (1 : 0 : 0)$. Parametrizing the curve using the non-singular point $(0:1:0)$ yields that
  \[d=-\frac{k^2(8 s - 41 t)^2 (64 s^2 - 272 s t + 97 t^2)^2}{(8 s - 
    25 t)^2 (64 s^2 - 784 s t + 2209 t^2)^2}, \quad z=\frac{k(8 s - 41 t) (8 s - 33 t) (64 s^2 - 272 s t + 97 t^2)}{8 (8 s - 
   25 t) t (64 s^2 - 784 s t + 2209 t^2)}, \quad s,t\in\Q\]
   where the orbit is 
 {\footnotesize  \[\frac{k(8 s - 41 t) (8 s - 33 t) (64 s^2 - 272 s t + 97 t^2)}{8 (8 s - 
   25 t) t (64 s^2 - 784 s t + 2209 t^2)}\longrightarrow \frac{k (8 s - 41 t)^2 (64 s^2 - 272 s t + 97 t^2)}{(8 s - 
   25 t)^2 (64 s^2 - 784 s t + 2209 t^2)}\longrightarrow -\frac{k (64 s^2 - 272 s t + 97 t^2)}{64 s^2 - 784 s t + 2209 t^2} \]}  
   {\footnotesize  \[ \longrightarrow -\frac{k (64 s^2 - 272 s t + 97 t^2)}{64 s^2 - 784 s t + 2209 t^2}.\]}  
The excluded value $t=0$ corresponds to $z=\infty$ and $d=-k^2$. In this case $\infty$ is indeed of type $1_2$, namely $\infty\rightarrow k\rightarrow -k\rightarrow -k$, and the other point of type $1_2$ is $0$, which is obtained for $8s=33t$. The fact that $\theta_{d,k}$ has exactly one $\Q$-rational fixed point follows from the observation that $-d$ is a $\Q$-rational square, whereas for $\theta_{d,k}$ to possess three $\Q$-rational fixed points, $-3d$ must be a $\Q$-rational square, Proposition \ref{prop:period1multiples}. 

   iii) The generalized $(3,1)$-dynatomic polynomial of $\theta_{d,k}$ is defined by  
   {\footnotesize\begin{equation*}
    \begin{split}
         \Phi^{*}_{\theta_{d,k},3,1}(z)&= -2(d-k^2)^4\big(d^8 + k^4 z^{12} + 6 d k^2 z^{10} (11 k^2 - 8 k z + z^2) +
  6 d^7 (k^2 - 8 k z + 11 z^2) \\ &+ 
  d^2 z^8 (495 k^4 - 880 k^3 z + 396 k^2 z^2 - 48 k z^3 + z^4) +  
  22 d^3 z^6 (42 k^4 - 144 k^3 z + 135 k^2 z^2 - 40 k z^3 + 3 z^4) \\ &+  
  99 d^4 z^4 (5 k^4 - 32 k^3 z + 56 k^2 z^2 - 32 k z^3 + 5 z^4) +  
  22 d^5 z^2 (3 k^4 - 40 k^3 z + 135 k^2 z^2 - 144 k z^3 + 42 z^4) \\ &+  
  d^6 (k^4 - 48 k^3 z + 396 k^2 z^2 - 880 k z^3 + 495 z^4)\big).
    \end{split}
    \end{equation*}}
    For $k\ne0$, applying the transformation $z\mapsto z/k$, $d\mapsto d/k^2$, the vanishing of $\Phi^{*}_{\theta_{d,k},3,1}(z)$ defines a planar curve in $\mathbb{P}^2$ with coordinates $(z:d:w)$. The curve is not geometrically irreducible over $\Q$. In fact, it consists of two irreducible components defined over the number field generated by a root of the polynomial $f(x)=x^2 + 4459x/32 + 19849913/4096$. One may show as in the proof of Theorem \ref{thm:morecycles} ii) that the only $\Q$-rational points on the curve are the singular points $(0 : 0 : 1), (1 : 1 : 1), (0 : 1 : 0), (1 : 0 : 0)$, none of which gives rise to a map $\theta_{d,k}$, $d\ne0$, $k^2\ne d$. For $k=0$, see Theorem \ref{thm:preper0} iii). Full computations can be found in \href{https://github.com/hasanbilgili44/Calculation-for-Paper/blob/main/Nonabelian-Automorphism/Theorem%206.2}{Theorem 6.2} of \cite{Github}.
\end{Proof}

\begin{Proposition}
\label{prop:preper0}
    Consider the map $\theta_{d,0}(z)=-2dz/(z^2+d)$ where  $d\in \Q\setminus \{0\}$.
    \begin{itemize}
        \item [i)] $\theta_{d,0}$ has no $\Q$-rational preperiodic point of type $2_m$, $m\ge 1$.
        \item[ii)] $\theta_{d,0}$ has no $\Q$-rational preperiodic point of type $3_m$, $m\ge1$.
        \item [iii)] For $n>3$, the map $\theta_{d,0}$ has no $\Q$-rational preperiodic points of type $n_m$, $m\ge 1$.
    \end{itemize}
\end{Proposition}
\begin{Proof}
   i) If $\theta_{d,0}$ has a $\Q$-rational periodic point of exact period $2$, then $d=c^2$ where $c$ and $-c$ are the two periodic points with period $2$, see Proposition \ref{prop:period2}. The result now follows from Remark \ref{rem:2cyclenotail}.
   
    ii) The generalized $(1,3)$-dynatomic polynomial of $\theta_{d,0}$ is given by 
    \begin{equation*}
     \Phi^*_{\theta_{d,0},1,3}= d^3 (d^3 + 33 d^2 z^2 + 27 d z^4 + 3 z^6).   
    \end{equation*}
    The statement follows from the irreducibility of $3Z^3+27Z^2+33Z+1$ over $\Q$.
    
    iii) The statement follows from Theorem \ref{thm}.
\end{Proof}

\begin{Theorem}
\label{thm:preper2}
Let $\theta_{d,k}(z)=(kz^2-2dz+dk)/(z^2-2kz+d)$ with $k\in \Q$, $d\in \Q\setminus \{0\}$, $k^2\ne d$. 
\begin{itemize}
\item[i)] $\theta_{d,k}$ has no $\Q$-rational preperiodic point of type $2_m$, $m\ge 1$.
\item[ii)] $\theta_{d,k}$, $k\ne 0$, has a $\Q$-rational preperiodic point $z$ of type $3_1$ if and only if {\footnotesize
 \[d=-\frac{k^2(512 s^3 - 10944 s^2 t + 59544 s t^2 - 17001 t^3)^2}{3 (512 s^3 - 17088 s^2 t + 183960 s t^2 - 640713 t^3)^2}, \quad z=\frac{k(8 s - 81 t) (512 s^3 - 10944 s^2 t + 59544 s t^2 - 17001 t^3)}{24 t (512 s^3 - 17088 s^2 t + 183960 s t^2 - 640713 t^3)},\quad t\ne 0.\]} In this case, $\theta_{d,k}$ possesses three $\Q$-rational preperiodic points of type $3_1$.
 \item[iii)] $\theta_{d,k}$ has no $\Q$-rational preperiodic points of type $3_m$, $m\ge 2$.
\item[iv)] For $n>3$, the map $\theta_{d,k}$ has no $\Q$-rational preperiodic points of type $n_m$, $m\ge 1$.
\end{itemize}
\end{Theorem}
\begin{Proof}
i) If $\theta_{d,k}$ has a $\Q$-rational periodic point of exact period $2$, then, according to Proposition \ref{prop:period2}, one has $d=c^2$ where $c$ and $-c$ are the two periodic points with period $2$. The result now follows from Remark \ref{rem:2cyclenotail}.

ii) The generalized $(1,3)$-dynatomic polynomial of $\theta_{d,k}$ is given by 

{\footnotesize\begin{eqnarray*}
\Phi^*_{\theta_{d,k},1,3}= (d-k^2)^2\big(d^4 + k^2 z^6 + 3 d z^4 (11 k^2 - 8 k z + z^2) + 
    3 d^3 (k^2 - 8 k z + 11 z^2) + 
    d^2 z^2 (27 k^2 - 80 k z + 27 z^2)\big).
\end{eqnarray*}}
When $k\ne0$, we may apply the transformation $z\mapsto z/k,d\mapsto d/k^2$ to obtain the planar curve 
$$d^4 + z^6 + 3 d z^4 (11 - 8 z + z^2) + 3 d^3 (1 - 8 z + 11 z^2) + 
 d^2 z^2 (27 - 80 z + 27 z^2)=0.$$
 The latter is a genus $0$ curve with three singular points $(d:z:w)=(0 : 0 : 1), (1 : 1 : 1), (0 : 1 : 0)$. Parametrizing the curve using the non-singular point $(1:0:0)$ yields that
{\footnotesize
 \[d=-\frac{k^2(512 s^3 - 10944 s^2 t + 59544 s t^2 - 17001 t^3)^2}{3 (512 s^3 - 17088 s^2 t + 183960 s t^2 - 640713 t^3)^2}, \quad z=\frac{k(8 s - 81 t) (512 s^3 - 10944 s^2 t + 59544 s t^2 - 17001 t^3)}{24 t (512 s^3 - 17088 s^2 t + 183960 s t^2 - 640713 t^3)},\quad t\ne 0.\]}
The excluded value $t=0$ corresponds to $z=\infty$ and $d=-k^2/3$; in this case $\infty$ is of type $3_1$, see the discussion following Theorem \ref{thm:period3}, and the other two points of type $3_1$ are finite. That $\theta_{d,k}$ possesses three distinct $\Q$-rational preperiodic points of type $3_1$ follows from Theorem \ref{thm:period3}. For $k=0$, there is no $\Q$-rational point of type $3_1$ since $\theta_{d,0}$ has no $\Q$-rational $3$-cycle.

iii) The generalized $(2,3)$-dynatomic polynomial of $\theta_{d,k}$ is given by 
{\footnotesize
\begin{equation*}
        \begin{split}
            \Phi^{*}_{\theta_{d,k},2,3}(z)
         &= (d - k^2)^4 \big(d^8 + k^4 z^{12} + 
   2 d k^2 z^{10} (69 k^2 - 48 k z + 7 z^2) + 
   2 d^7 (7 k^2 - 48 k z + 69 z^2)  
  \\&+  d^2 z^8 (975 k^4 - 1760 k^3 z + 780 k^2 z^2 - 96 k z^3 + z^4) 
  \\&+ 2 d^3 z^6 (934 k^4 - 3168 k^3 z + 2985 k^2 z^2 - 880 k z^3 + 69 z^4)
  \\&+ 2 d^5 z^2 (69 k^4 - 880 k^3 z + 2985 k^2 z^2 - 3168 k z^3 + 934 z^4) 
  \\&+ d^4 z^4 (975 k^4 - 6336 k^3 z + 11048 k^2 z^2 - 6336 k z^3 + 975 z^4) 
  \\&+ d^6 (k^4 - 96 k^3 z + 780 k^2 z^2 - 1760 k z^3 + 975 z^4)\big).
        \end{split}
    \end{equation*}}
 For $k\ne0$, applying the transformation $z\mapsto z/k$, $d\mapsto d/k^2$, the vanishing of $\Phi^{*}_{\theta_{d,k},2,3}(z)$ defines a planar curve in $\mathbb{P}^2$ with coordinates $(z:d:w)$. The curve is not geometrically irreducible over $\Q$. It consists of two irreducible components defined over the number field generated by a root of the  polynomial $f(x)=x^2 + 4715x/32 + 22034617/4096$. One may show, as in the proof of Theorem \ref{thm:morecycles} ii), that the only $\Q$-rational points on the curve are the singular points $(0 : 0 : 1), (1 : 1 : 1), (0 : 1 : 0), (1 : 0 : 0)$. The latter points do not give rise to a map $\theta_{d,k}$, $d\ne0$, $k^2\ne d$. For $k=0$, there is no $\Q$-rational $3$-cycle. For full details of the latter computations, see \href{https://github.com/hasanbilgili44/Calculation-for-Paper/blob/main/Nonabelian-Automorphism/Theorem%206.4}{Theorem 6.4} of \cite{Github}.
 
iv) The statement follows from Theorem \ref{thm}.
\end{Proof}

\begin{Corollary}
    
 \label{thm : preper 1_1 and period} 
    Let $\theta_{d,k}(z)=(kz^2-2dz+dk)/(z^2-2kz+d)$ with $k\in \Q$, $d\in \Q\setminus \{0\}$, $k^2\ne d$. Assume that $k \ne 0$. $\theta_{d,k}$ possesses both a $\Q$-rational preperiodic point of type $1_1$ together with a $\Q$-rational periodic point of exact period $2$ if and only if $d=c^2$, $k=c^2(c^2+3r^2)/(r(3c^2+r^2))$, $r,c\in\Q\setminus\{0\}$. In this case, one has $\theta_{d,k}(r)=c^2/r$, $\theta_{d,k}(c^2/r)=c^2/r$; and $\theta_{d,k}(\pm c)=\mp c$. In addition, $\theta_{d,k}$ possesses exactly one $\Q$-rational fixed point.
\end{Corollary}
\begin{Proof}
    By Theorem \ref{thm:preper1} i), a $\Q$-rational point of type $1_1$ exists if and only if $\theta_{d,k}$ has a $\Q$-rational fixed point $v$, and $v\ne 0$ since $k\ne0$. By Proposition \ref{prop:period2}, $d=c^2$. Solving $c^2=-(3k-v)v^2/(k-3v)$, see Proposition \ref{prop:period1}, for $k$ gives $k=v(v^2+3c^2)/(c^2+3v^2)$, which is the stated formula with $r=c^2/v$. The fixed point is unique by Theorem \ref{thm:morecycles} iii).
\end{Proof}

\begin{Corollary}
    
\label{thm : preper 1_2 and period} 
    Let $\theta_{d,k}(z)=(kz^2-2dz+dk)/(z^2-2kz+d)$ with $k\in \Q$, $d\in \Q\setminus \{0\}$, $k^2\ne d$. Assume that $k \ne 0$.
    \begin{itemize}
        \item [i)] If $\theta_{d,k}$ possesses a $\Q$-rational preperiodic point of type $1_2$, then it does not possess a $\Q$-rational periodic point of exact period $2$. 
        \item[ii)] If $\theta_{d,k}$ possesses a $\Q$-rational preperiodic point of type $1_2$, then it does not possess a $\Q$-rational periodic point of exact period $3$. 
    \end{itemize}
\end{Corollary} 
\begin{Proof}
    We observe that the existence of a $\Q$-rational preperiodic point of type $1_2$ forces $-d$ to be a $\Q$-rational square, see Theorem \ref{thm:preper1} ii); whereas the existence of a $\Q$-rational periodic point of exact period $2$ (exact period $3$, respectively) implies that $d$ must be a $\Q$-rational square ($-3d$ must be a $\Q$-rational square, respectively), see Proposition \ref{prop:period2} (see Theorem \ref{thm:period3}, respectively). This concludes the proof.
\end{Proof}

We are now in a position to classify the portraits $\mathcal P(\theta_{d,k},\Q)$.

\begin{Theorem}
\label{thm:classification}
Let $\theta_{d,k}(z)=(kz^2-2dz+dk)/(z^2-2kz+d)$ with $k\in \Q$, $d\in \Q\setminus \{0\}$, $k^2\ne d$. Then $\mathcal P(\theta_{d,k},\Q)$ is one of the seven portraits in Table \ref{Table1}. More precisely, according to the set of $\Q$-rational cycles of $\theta_{d,k}$, one of the following holds.
\begin{itemize}
\item[a)] If $\theta_{d,k}$ has no $\Q$-rational periodic point, then $\Preper(\theta_{d,k},\Q)=\emptyset$, and $\mathcal P(\theta_{d,k},\Q)$ is the portrait vii.
\item[b)] If the $\Q$-rational periodic points of $\theta_{d,k}$ consist of exactly one fixed point, then $\mathcal P(\theta_{d,k},\Q)$ is the portrait i or iv.
\item[c)] If $\theta_{d,k}$ has three $\Q$-rational fixed points, then $\mathcal P(\theta_{d,k},\Q)$ is the portrait iii.
\item[d)] If the $\Q$-rational periodic points of $\theta_{d,k}$ consist of exactly one $2$-cycle, then $\mathcal P(\theta_{d,k},\Q)$ is the portrait vi.
\item[e)] If the $\Q$-rational periodic points of $\theta_{d,k}$ consist of a fixed point together with a $2$-cycle, then $\mathcal P(\theta_{d,k},\Q)$ is the portrait ii.
\item[f)] If $\theta_{d,k}$ has a $\Q$-rational periodic point of exact period $3$, then $\mathcal P(\theta_{d,k},\Q)$ is the portrait v.
\end{itemize}
\end{Theorem}
\begin{Proof}
By Theorem \ref{thm}, every $\Q$-rational periodic point of $\theta_{d,k}$ has exact period $1$, $2$ or $3$. Since $\theta_{d,k}(\infty)=k\ne\infty$, the fixed points of $\theta_{d,k}$ are the roots of the cubic polynomial (\ref{eq1}), whose discriminant is $-108d(d-k^2)^2\ne 0$. Since its roots are distinct, the number of $\Q$-rational fixed points of $\theta_{d,k}$ is $0$, $1$ or $3$. By Theorem \ref{thm:morecycles} ii)--v), Remark \ref{rem:morecycles0}, and Theorem \ref{thm:period3}, the set of $\Q$-rational cycles of $\theta_{d,k}$ is exactly as described in one of the cases a)--f). Since every $\Q$-rational preperiodic point is mapped after finitely many iterations into a $\Q$-rational cycle, it remains to determine the $\Q$-rational points of type $n_m$, $m\ge1$, in each case. We will use repeatedly that if there are no $\Q$-rational points of type $n_m$, then there are no $\Q$-rational points of type $n_{m+1}$. Moreover, by Lemma \ref{lem:sameimage} i), the preimages of any $w\in\PP^1(\Q)$ are of the form $z$ and $d/z$; so if one preimage is $\Q$-rational, then so is the other one.

a) This is immediate.

b) Let $z$ be the fixed point. By Theorem \ref{thm:preper1} i), $d/z$ is the unique $\Q$-rational point of type $1_1$. The $\Q$-rational points of type $1_2$ are the preimages $w$, $d/w$ of $d/z$. If they exist, then $-d$ is a $\Q$-rational square, see Theorem \ref{thm:preper0} ii) and Theorem \ref{thm:preper1} ii); in particular $d$ is not a $\Q$-rational square, hence $w\ne d/w$. Points of type $1_m$, $m\ge 3$, do not exist by Theorem \ref{thm:preper0} iii) and Theorem \ref{thm:preper1} iii). Therefore, $\mathcal P(\theta_{d,k},\Q)$ is the portrait i if there is no $\Q$-rational point of type $1_2$, and the portrait iv otherwise.

c) Let $z_1,z_2,z_3$ be the fixed points. By Theorem \ref{thm:preper1} i), $d/z_i$ is a $\Q$-rational point of type $1_1$, $i=1,2,3$, and these three points are pairwise distinct. If $k\ne 0$, then there are no $\Q$-rational points of type $1_2$ by Theorem \ref{thm:preper1} ii). If $k=0$, then $-3d$ is a $\Q$-rational square, see Proposition \ref{prop:period1multiples}, so $-d$ is not a $\Q$-rational square, and Theorem \ref{thm:preper0} ii) implies that there are no $\Q$-rational points of type $1_2$. Hence, $\mathcal P(\theta_{d,k},\Q)$ is the portrait iii.

d) By Proposition \ref{prop:period2}, $d=c^2$ and the $2$-cycle is $\{c,-c\}$. According to Remark \ref{rem:2cyclenotail}, there are no $\Q$-rational points of type $2_m$, $m\ge1$. Thus, $\mathcal P(\theta_{d,k},\Q)$ is the portrait vi.

e) Let $z$ be the fixed point. As in d), the $2$-cycle has no $\Q$-rational tail. By Theorem \ref{thm:preper1} i), $d/z$ is the unique $\Q$-rational point of type $1_1$. If $k\ne 0$, then there are no $\Q$-rational points of type $1_2$ by Corollary \ref{thm : preper 1_2 and period} i). If $k=0$, then $d$ is a $\Q$-rational square by Proposition \ref{prop:period2}, so $-d$ is not, and there are no $\Q$-rational points of type $1_2$ by Theorem \ref{thm:preper0} ii). Hence, $\mathcal P(\theta_{d,k},\Q)$ is the portrait ii.

f) By Theorem \ref{thm:period3}, $k\ne 0$, $-3d$ is a $\Q$-rational square, and for every point $z$ in the $3$-cycle, $d/z$ is a $\Q$-rational point of type $3_1$. These three points are pairwise distinct. Moreover, none of them lies in the $3$-cycle: if $d/z_i=z_j$ for points $z_i,z_j$ in the $3$-cycle, then $\theta_{d,k}(z_i)=\theta_{d,k}(z_j)$, so $z_i=z_j=d/z_i$. This is impossible if $z_i\in\{0,\infty\}$, while for $z_i\in\Q\setminus\{0\}$ it gives $d=z_i^2$, which is impossible since $-3d$ is a $\Q$-rational square. (The $3$-cycle may contain $\infty$, see the discussion following Theorem \ref{thm:period3}; in this case the corresponding point of type $3_1$ is $d/\infty=0$.) By Theorem \ref{thm:preper2} iii), there are no $\Q$-rational points of type $3_m$, $m\ge 2$. Therefore, $\mathcal P(\theta_{d,k},\Q)$ is the portrait v.
\end{Proof}

\begin{ProofOf}{Theorem \ref{thm:conj}}
By Lemma \ref{lem:Yasuf}, $f$ is $\Q$-conjugate to a map $\theta_{d,k}$, and $\Q$-conjugate maps have isomorphic portraits. By Theorem \ref{thm:classification}, $\mathcal P(\theta_{d,k},\Q)$ is one of the portraits in Table \ref{Table1}, which have $2,4,6,4,6,2,0$ vertices, respectively. Hence $|\Preper(f,\Q)|\le 6$. The bound is attained, see Example \ref{ex:portraits} iii) and v).
\end{ProofOf}

\begin{Remark}
\label{rem:bare2cycle}
Case d) of Theorem \ref{thm:classification} occurs, see Example \ref{ex:portraits} vi). More precisely, in view of Proposition \ref{prop:period2}, Proposition \ref{prop:period1}, Theorem \ref{thm:morecycles} v) and Corollary \ref{thm : preper 1_1 and period}, the map $\theta_{c^2,k}$, $c\in\Q\setminus\{0\}$, $k\ne\pm c$, realizes the portrait vi if and only if $k\ne 0$ and $k\ne c^2(c^2+3r^2)/(r(3c^2+r^2))$ for every $r\in\Q\setminus\{0\}$. Indeed, $0$ is a fixed point of $\theta_{c^2,k}$ if and only if $k=0$, and a nonzero fixed point $v\in\Q$ corresponds to $r=c^2/v$ in Corollary \ref{thm : preper 1_1 and period}.
\end{Remark}

We now produce examples of maps $\theta_{d,k}$ with every possible portrait associated with $\Preper(\theta_{d,k},\Q)$, $k\in \Q$, $d\in \Q\setminus \{0\}$, $k^2\ne d$, as in Table \ref{Table1}.   
\begin{Example}
\label{ex:portraits}
In each of the following examples, we display the portrait associated with $\Preper(\theta_{d,k},\Q)$ for the given  pair $(d,k)$.  
\begin{itemize} 
\item[i)] $(d,k)=\left(2,\frac{7}{5}\right)$
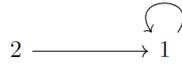
\begin{center}
\begin{tikzcd}
    2 \arrow[r] & 1 \arrow[out=-45, in=45, loop, distance=2em]
\end{tikzcd}
\end{center}
\item[ii)] $(d,k)= \left(\frac{49}{324},-\frac{1}{2}\right)$
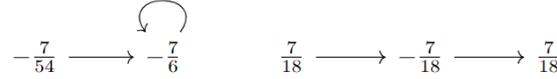
\begin{center}
\begin{tikzcd}
    -\frac{7}{54} \arrow[r] & -\frac{7}{6} \arrow[out=-45, in=45, loop, distance=2em]
\end{tikzcd}
\hspace{1cm}
\begin{tikzcd}
    \frac{7}{18} \arrow[r, bend left] & -\frac{7}{18} \arrow[l, bend left]
\end{tikzcd}
\end{center}
\item[iii)] $(d,k)=\left(-\frac{100}{2187},-\frac{1}{3}\right)$
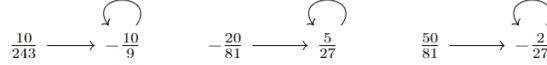
\begin{center}
\begin{tikzcd}
    \frac{10}{243} \arrow[r] & -\frac{10}{9} \arrow[out=-45, in=45, loop, distance=2em]
\end{tikzcd}
\hspace{0.6cm}
\begin{tikzcd}
    -\frac{20}{81} \arrow[r] & \frac{5}{27} \arrow[out=-45, in=45, loop, distance=2em]
\end{tikzcd}
\hspace{0.6cm}
\begin{tikzcd}
    \frac{50}{81} \arrow[r] & -\frac{2}{27} \arrow[out=-45, in=45, loop, distance=2em]
\end{tikzcd}
\end{center}
$(d,k)=\left(-3a^2,0\right)$, $a\in \Q\setminus\{0\}$
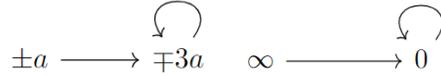
\begin{center}
\begin{tikzcd}
    a \arrow[r] & -3a \arrow[out=-45, in=45, loop, distance=2em]
\end{tikzcd}
\hspace{0.6cm}
\begin{tikzcd}
    -a \arrow[r] & 3a \arrow[out=-45, in=45, loop, distance=2em]
\end{tikzcd}
\hspace{0.6cm}
\begin{tikzcd}
    \infty \arrow[r] & 0 \arrow[out=-45, in=45, loop, distance=2em]
\end{tikzcd}
\end{center}
\item[iv)] $(d,k)=\left(- \frac{13417569}{640747969},1\right)$
\begin{center}
\begin{tikzcd}[row sep=small]
    \frac{91575}{202504} \arrow[r] & -\frac{120879}{430321} \arrow[r] & \frac{111}{1489} \arrow[out=-45, in=45, loop, distance=2em] \\
    & -\frac{29304}{632825} \arrow[u] &
\end{tikzcd}
\end{center}
$(d,k)=(-a^2,0)$, $a\in \Q\setminus\{0\}$
\begin{center}
\begin{tikzcd}[row sep=small]
    a \arrow[r] & \infty \arrow[r] & 0 \arrow[out=-45, in=45, loop, distance=2em] \\
    & -a \arrow[u] &
\end{tikzcd}
\end{center}
\item[v)] $(d,k)=\left(-\frac{1031116321}{672121026723},1\right)$
\begin{center}
\begin{tikzcd}[row sep=large, column sep=large]
    & & \frac{2087215}{45912913} \arrow[dd] & \frac{1573439}{115018947} \arrow[l] \\
    \frac{2344103}{11359896} \arrow[r] & -\frac{2600991}{23193121} \arrow[ur, bend left=25] & & \\
    & & -\frac{256888}{34553017} \arrow[ul, bend left=25] & -\frac{3114767}{92299155} \arrow[l]
\end{tikzcd}
\end{center}
\item[vi)] $(d,k)=(1,2)$
\begin{center}
\begin{tikzcd}
    1 \arrow[r, bend left] & -1 \arrow[l, bend left]
\end{tikzcd}
\end{center}
Indeed, $\theta_{1,2}(1)=-1$ and $\theta_{1,2}(-1)=1$. The polynomial $2-3z+6z^2-z^3$ in (\ref{eq1}) has no rational root, and $-3d=-3$ is not a $\Q$-rational square, so $\theta_{1,2}$ has no $\Q$-rational fixed point and no $\Q$-rational periodic point of exact period $3$, see Theorem \ref{thm:period3}. The claim now follows from Theorem \ref{thm:classification} d).
\item[vii)] $(d,k)=(2,1)$
\begin{center}
$\Preper(\theta_{2,1},\Q)=\emptyset$
\end{center}
Indeed, the polynomial $2-6z+3z^2-z^3$ in (\ref{eq1}) has no rational root, $d=2$ and $-3d=-6$ are not $\Q$-rational squares, so $\theta_{2,1}$ has no $\Q$-rational periodic points, see Proposition \ref{prop:period2}, Theorem \ref{thm:period3} and Theorem \ref{thm}.
\end{itemize}
\end{Example}

\section*{Conflict of Interest statement} On behalf of all authors, the corresponding author states that there is no conflict of interest.

\section*{Data availability statement} The authors declare that the data supporting the findings of this study are available within the paper and in the GitHub repository \cite{Github}.

\end{document}